%% file: agt-5-1.tex
\documentclass{gtart_h}
\input agtout

\lognumber{1}
\volumenumber{5}
\volumeyear{2005}
\papernumber{1}
\pagenumbers{1}{22}
\received{13 October 2004} 
\revised{6 November 2004}
\accepted{7 December 2004}
\published{5 January 2005}

\usepackage{amsmath,amssymb,graphicx,psfrag}

\allowdisplaybreaks

\theoremstyle{plain}
 \newtheorem{thm}{Theorem}[section]
 \newtheorem{prop}[thm]{Proposition}
\newtheorem{lem}[thm]{Lemma}
 \newtheorem{cor}[thm]{Corollary}
 
\theoremstyle{definition}
 \newtheorem{defn}{Definition}
 \newtheorem{rmk}{Remark}
 \newtheorem{example}{Example}

\newcommand{\Z}{\mathbb{Z}}
\newcommand{\C}{\mathbb{C}}
\newcommand{\R}{\mathbb{R}}
\newcommand{\s}{\sigma}
\renewcommand{\a}{\mathbf{a}}
\newcommand{\x}{\mathbf{x}}
\renewcommand{\c}{\cite} 
\newcommand{\tZ}{\mathbb{Q}(A)}
\newcommand{\ZA}{{\mathbb{Z}}[A^{\pm1}] }
\newcommand{\SL}{\mathfrak s\mathfrak l}

\newcommand{\kb}[1]{\ensuremath{\langle #1 \rangle}}

\newcommand{\myone}{\ensuremath{{\bf 1}}}
\newcommand{\tr}{{\rm tr}}
\renewcommand\theenumi{{\it \roman{enumi}}} 
\newcommand{\splice}{\ensuremath\kb {\makebox[0.3cm][c]{\raisebox{-0.3ex}{\rotatebox{90}{$\asymp$}}}}}

\begin{document}

\title{On the Mahler measure of\\ Jones polynomials under twisting}
\authors{Abhijit Champanerkar\\Ilya Kofman}
\address{Department of Mathematics, Barnard College, Columbia University\\
3009 Broadway, New York, NY 10027, USA}
\secondaddress{Department of Mathematics, Columbia University\\
2990 Broadway, New York, NY 10027, USA}

\asciiaddress{Department of Mathematics, Barnard College, 
Columbia University\\3009 Broadway, New York, 
NY 10027, USA\\and\\Department of Mathematics, 
Columbia University\\2990 Broadway, New York, NY 10027, USA}

\gtemail{\mailto{abhijit@math.columbia.edu}, 
\mailto{ikofman@math.columbia.edu}}
\asciiemail{abhijit@math.columbia.edu, ikofman@math.columbia.edu}

\begin{abstract}
We show that the Mahler measures of the Jones polynomial and of the
colored Jones polynomials converge under twisting for any link.
Moreover, almost all of the roots of these polynomials approach the
unit circle under twisting.  In terms of Mahler measure convergence,
the Jones polynomial behaves like hyperbolic volume under Dehn
surgery.  For pretzel links $\mathcal{P}(a_1,\ldots,a_n)$, we show
that the Mahler measure of the Jones polynomial converges if all
$a_i\to\infty$, and approaches infinity for $a_i=$ constant if
$n\to\infty$, just as hyperbolic volume.  We also show that after
sufficiently many twists, the coefficient vector of the Jones
polynomial and of any colored Jones polynomial decomposes into fixed
blocks according to the number of strands twisted.
\end{abstract}

\asciiabstract{%
We show that the Mahler measures of the Jones polynomial and of the
colored Jones polynomials converge under twisting for any link.
Moreover, almost all of the roots of these polynomials approach the
unit circle under twisting.  In terms of Mahler measure convergence,
the Jones polynomial behaves like hyperbolic volume under Dehn
surgery.  For pretzel links P(a_1,...,a_n), we show that the Mahler
measure of the Jones polynomial converges if all a_i tend to infinity,
and approaches infinity for a_i = constant if n tend to infinity, just
as hyperbolic volume.  We also show that after sufficiently many
twists, the coefficient vector of the Jones polynomial and of any
colored Jones polynomial decomposes into fixed blocks according to the
number of strands twisted.}

\primaryclass{57M25}
\secondaryclass{26C10}
\keywords{Jones polynomial, Mahler measure, Temperley-Lieb algebra, hyperbolic volume}
\maketitle

\section{Introduction}

It is not known whether any natural measure of complexity of
Jones-type polynomial invariants of a knot is related to the volume of
the knot complement, a measure of the knot's geometric complexity.
The Mahler measure, which is the geometric mean on the unit circle, is in a sense the canonical measure of complexity on the space of polynomials \c{Skoruppa}: 
For any monic polynomial $f$ of degree $n$, let $f_k$ denote the polynomial whose roots are $k$-th powers of the roots of $f$, then for {\em any} norm $||\cdot ||$ on the vector space of degree $n$ polynomials,
$$ \lim_{k\to\infty} ||f_k||^{1/k} = M(f) $$ 
In this work, we show that the Mahler measure of the Jones polynomial and of the colored Jones polynomials behaves like hyperbolic volume under Dehn surgery.
Suppose $L_m$ is obtained from the hyperbolic link $L$ by adding $m$ full twists on $n$ strands of $L$.  
In other words, there is an unknot $U$ encircling $n$ strands of $L$, such that $L_m$ is obtained from $L$ by a $-1/m$ surgery on $U$.
By Thurston's hyperbolic Dehn surgery theorem,
\[ \lim_{m\to\infty} {\rm Vol}(S^3\setminus L_m) = {\rm Vol}(S^3\setminus (L\cup U)) \]
We show that as $m\to\infty$, the Mahler measure of the Jones polynomial of $L_m$ converges to the Mahler measure of a $2$-variable polynomial (Theorem \ref{thm1}).
We also show that for any $N$, the Mahler measure of the colored Jones polynomials $J_N({L_m};t)$ converges as $m\to\infty$ (Theorem \ref{colJPthm}).  
Moreover, as $m\to\infty$, almost all of the roots of these polynomials approach the unit circle (Theorem \ref{zeros}).
This result explains experimental observations by many authors who have studied the distribution of roots of Jones polynomials for various families of knots and links (eg, \c{Z2,Z1,Z4,Z3}).
All of our results extend to multiple twisting: adding $m_i$ twists on $n_i$ strands of $L$, for $i=1,\ldots,k$ (eg, Corollary \ref{multitwist}).

For pretzel links $\mathcal{P}(a_1,\ldots,a_n)$, we show that the Mahler measure of the Jones polynomial, 
$M(V_{\mathcal{P}(a_1,\ldots,a_n)}(t)),$ converges in one parameter and
approaches infinity in the other, just like hyperbolic volume:
\begin{itemize}
\item For odd $a_1= \ldots =a_n, \;
  M(V_{\mathcal{P}(a_1,\ldots,a_n)}(t)) \to \infty$ as $n \to \infty$
  (Theorem \ref{pretzelthm}).  Volume $\to\infty$ by a result of
  Lackenby \c{Lackenby}.
\item If all $a_i\to \infty, \; M(V_{\mathcal{P}(a_1,\ldots,a_n)}(t))$
  converges if $n$ is fixed (Corollary \ref{multitwist}).  Volume also
  converges by Thurston's Dehn surgery theorem.
\end{itemize}

In addition, we show that after sufficiently many twists on $n$
strands, the coefficient vector of the Jones polynomial decomposes
into $([n/2]+1)$ fixed blocks separated by zeros if $n$ is odd, and by
alternating constants if $n$ is even (Theorem \ref{blocks}).  For
example, these are the coefficient vectors of the Jones polynomial of
a knot after twisting the same $5$ strands (see Table \ref{t5}):

{\small{\bf $5$ full twists:}

\framebox{1  -1   2  -1   2  -1   1}   0   0   0   0   0   0   0   0   0   0   0   0   0   0   0   0   0   \framebox{2  -8  16 -23  20 -12   0   7  -7   4  -1}   0   0   0   \framebox{1  -5  15 -29  40 -42  33 -19   8  -2}

{\bf $20$ full twists:}

\framebox{1  -1   2  -1   2  -1   1}   0   0   0   0   0   0   0   0   0   0   0   0   0   0   0   0   0   0   0   0   0   0   0   0   0   0   0   0   0   0   0   0   0   0   0   0   0   0   0   0   0   0   0   0   0   0   0   0   0   0   0   0   0   0   0   0   0   0   0   0   0   0   0   0   0   0   0   0   0   0   0   0   0   0   0   0   0   0   0   0   0   0   0   0   0   0   0   0   0   0   0   0   \framebox{2  -8  16 -23  20 -12   0   7  -7   4  -1}   0   0   0   0   0   0   0   0   0   0   0   0   0   0   0   0   0   0   0   0   0   0   0   0   0   0   0   0   0   0   0   0   0   0   0   0   0   0   0   0   0   0   0   0   0   0   0   0   \framebox{1  -5  15 -29  40 -42  33 -19   8  -2}}

These fixed blocks completely determine the Jones polynomial under
twisting, simply moving apart linearly in the number of twists.
Moreover, this result extends to the colored Jones polynomial: If we
fix both $N$ and the number of strands twisted, then for $m$
sufficiently large, the coefficient vector of $J_N({L_m};t)$
decomposes into fixed blocks separated by zeros, which move apart
linearly with more twists (Corollary \ref{colblocks}).

This work is motivated by a deep open problem: how to bridge the chasm
between quantum and geometric topology.  Thurston established the
importance of geometric invariants, especially hyperbolic volume, in
low-dimensional topology.  Yet the vast families of quantum
invariants, which followed the discovery 20 years ago of the Jones
polynomial, are not understood in terms of geometry.  The ``volume
conjecture'' and its variants propose that colored Jones polynomials,
which are weighted sums of Jones polynomials of cablings, determine
the volume of hyperbolic knots (see \c{Gukov, MM}).  The original
Jones polynomial of a knot is still not understood in terms of the
knot complement.  A direct link between the Jones polynomial and the
volume of the knot complement would relate the most important quantum
and geometric invariants.  With that goal in mind, we compare the
volume, a measure of geometric complexity of the knot complement, with
the Mahler measure of the Jones polynomial, a natural measure of
complexity on the space of polynomials.
 
Experimental evidence on hyperbolic knots with a simple hyperbolic
structure is suggestive.  In \c{CKP}, we computed Jones polynomials
for hyperbolic knots whose complements can be triangulated with seven
or fewer ideal tetrahedra.  A glance at these polynomials reveals how
different they are from Jones polynomials in knot tables organized by
crossing number.  The span of the Jones polynomial gives a lower bound
for the crossing number, with equality for alternating knots.  The
spans of these polynomials vary from $4$ to $43$, but the polynomials
with large span are very sparse, and their nonzero coefficients are
very small.  Mahler measure is a natural measure on the space of
polynomials for which these kinds of polynomials are simplest.

The Mahler measure of other knot polynomials has been related to the
volume of the knot complement.  Boyd and Rodriguez-Villegas found
examples of knot complements (and other 1-cusped hyperbolic manifolds)
such that Vol$(M)=\pi\,{\rm m}(A)$, where $A(x,y)$ is the A-polynomial
\c{BRV,brv03}.  Silver and Williams \c{SW1} showed the Mahler measure
of the Alexander polynomial converges under twisting just like the
volume converges under the corresponding Dehn surgery: If $U$ has
nonzero linking number with some component of $L$, and $L_m$ is
obtained from $L$ by a $-1/m$ surgery on $U$, then the Mahler measure
of the multivariable Alexander polynomial of $L_m$ converges to that
of $L\cup U$.  In contrast, for the Jones polynomial, the limit in
Theorem \ref{thm1} is not the Mahler measure of the Jones polynomial
of $L\cup U$, according to all our examples.

This paper is organized as follows.
In Section \ref{sec2}, we show that the Mahler measure of Jones polynomials and colored Jones polynomials converges under twisting and multi-twisting, using the representation theory of braid groups and linear skein theory.  As a consequence, we show that almost all of the roots of these polynomials approach the unit circle.  
In Section \ref{twistsection}, we show that the coefficient vector of these polynomials decomposes after sufficiently many twists.
In Section \ref{pretzelsection}, we show that for a family of pretzel links, the Mahler measure of the Jones polynomial both converges and diverges like volume.
We conclude with some observations using the Knotscape census of knots up to $16$ crossings.

\subsection*{Acknowledgments}
We thank Vaughan Jones and Adam Sikora for very helpful discussions on
related aspects of representation theory.  We thank Walter Neumann and
Andrzej Schinzel for suggesting Lemma 1 on page 187 of \c{Schinzel2}
used in the proof of Theorem \ref{zeros}.  We also thank David Boyd
for the idea of the proof of Lemma \ref{Boyd}.

\section{Mahler measure convergence}\label{sec2}

\begin{defn}\label{MahlerDef}
Let $f \in \C[z^{\pm 1}_1,\ldots,z^{\pm 1}_s]$.  The Mahler measure of $f$ is defined as follows, where $\exp(-\infty)=0$,
\[ M(f) = \exp\int_0^1\cdots\int_0^1\log\left|f(e^{2\pi i\theta_1},\ldots,e^{2\pi i\theta_s})\right| {\rm d} \theta_1\cdots {\rm d} \theta_s \]
The logarithmic Mahler measure is ${\rm m}(f)=\log M(f)$.
\end{defn}
The Mahler measure is multiplicative, $M(f_1 f_2) = M(f_1) M(f_2)$, and the logarithmic Mahler measure is additive.
If $s=1,\; f(z)=a_0 z^k\prod_{i=1}^n(z-\alpha_i)$, then by Jensen's formula, 
\begin{equation}\label{Jensen}
M(f) = |a_0|\prod_{i=1}^n {\rm max}(1,|\alpha_i|)
\end{equation}
For vectors $\a,\x\in\Z^s$, let $h(\x)=\max|x_i|$ and 
\[ \nu(\a) = \min\{ h(\x)\;|\; \x\in\Z^s,\; \a\cdot\x=0 \} \]
For example, $\nu(1,d,\ldots,d^{s-1})=d$.
\begin{lem}[Boyd, Lawton \c{Lawton}]\label{BL}
For every $f \in \C[z^{\pm 1}_1,\ldots,z^{\pm 1}_s]$,
\[ M(f) = \lim_{\nu(\x)\to\infty} M\left(f(z^{x_1},\ldots,z^{x_s})\right) \]
with the following useful special case:
$M(f) = \lim\limits_{d\to\infty} M\left(f(z,z^d,\ldots,z^{d^{s-1}})\right)$.
\end{lem}
For a good survey on the Mahler measure of polynomials, see \c{Schinzel}.
Finally, because the Mahler measure is multiplicative, Definition \ref{MahlerDef} can be naturally extended to rational functions of Laurent polynomials.

The Temperley-Lieb algebra $TL_n$ is closely related to the Jones polynomial.
$TL_n$ is the algebra over $\ZA$ with generators $\{\myone, e_1, e_2,\ldots,e_{n-1}\}$ and relations, with $\delta=-A^2-A^{-2}$,
\begin{equation}\label{TL}
e_i^2=\delta e_i,\quad e_i e_{i\pm 1} e_i = e_i,\quad e_i e_j = e_j e_i \text{ if } |i-j|\geq 2
\end{equation}
Kauffman gave a diagramatic interpretation of the Jones representation of the braid group, $\rho: B_n \to TL_n$ by $\rho(\sigma_i)=A\myone+A^{-1}e_i$, with the Markov trace interpreted as the bracket polynomial of the closed braid \c{Ktams}:
\[\tr(\rho(\beta)) = \langle\bar{\beta}\rangle  \]
Generalizing from braids to tangles using the Kauffman bracket skein relations, $TL_n$ is precisely the skein algebra of $D^2$ with $2n$ marked points on the boundary, with coefficients in $\ZA$.
The basis as a free $\ZA$-module consists of all diagrams with no crossings and no closed curves. 
Its dimension is the Catalan number $C_n = \frac{1}{n+1} \binom{2n}{n}$.
If the disc is considered as a square with $n$ marked points on the left edge and $n$ on the right, the product in the algebra is given by juxtaposing two squares to match marked points on the left edge of one square with the marked points on the right edge of the other square.

Moreover, since the skein algebras of $\R^2$ and $S^2$ are naturally isomorphic, a bilinear pairing is induced from the decomposition of $S^2$ into complementary discs $D'\cup D''$.
For any link diagram $L$ in $S^2$, decompose $L=L'\cup L''$ such that $L'=L\cap D',\, L''=L\cap D''$ and $L$ intersects the boundary of the disc in $2n$ points away from the crossings.
The bilinear pairing 
\[\langle\ ,\,\rangle : TL_n \times TL_n \to \ZA\] 
is given by $\langle L',L''\rangle  = \langle L\rangle $.
For a detailed introduction, see \c{Prasolov-Sossinsky}.

At this point, it is useful to return to Jones' orignial construction, in which $TL_n$ is viewed as a quotient of the Hecke algebra.
In \c{J}, Jones showed that irreducible representations of $B_n$ are indexed by Young diagrams, and the ones with at most two columns are the $TL_n$ representations.
By abuse of notation, let $TL_n$ also denote the algebra with coefficients in $\tZ$, and generators and relations as in (\ref{TL}).
The $\tZ$-basis is the same $\ZA$-basis described above, so the bilinear pairing extends linearly to $\tZ$ such that for any tangles, $\langle L',L''\rangle\in\ZA$.
With the ground ring $\ZA$ extended to $\tZ$,
$TL_n$ is semisimple, which implies that every irreducible submodule is generated by a minimal idempotent, and these idempotents can be chosen to be mutually orthogonal.  
More explicitly, Wenzl constructed orthogonal representations of Hecke algebras and gave an inductive formula for minimal idempotents such that $p_Y p_{Y'} = 0$ if the Young diagrams $Y\neq Y'$, and $\sum p_Y = 1$ (Corollary 2.3 \c{Wenzl}).

Let $L$ be any link diagram.  We deform $L$ to a link $L_m$ by performing a $-1/m$ surgery on an unknot $U$ which encircles $n$ strands of $L$.  This is the same as adding $m$ full right twists on $n$ strands of $L$.  It is useful to restate this in terms of the decomposition above, where $L \subset S^2 = D'\cup D''$.
The full right twist on $n$ strands in $B_n$ is denoted by $\Delta^2 = (\s_1\ldots\s_{n-1})^n$.
Let $L=\mathbf{1}_n\cup L''$ be any link diagram such that $\mathbf{1}_n$ is just the trivial braid on $n$ strands, and $\mathbf{1}_n=L\cap D',\, L''=L\cap D''$.
Let $L_m$ be the link obtained from $L$ by changing $\mathbf{1}_n$ to be $m$ full twists $\Delta^{2m}$, and leaving $L''$ unchanged.

\begin{thm}\label{thm1}
The Mahler measure of the Jones polynomial of $L_m$ converges as $m\to\infty$ to the Mahler measure of a $2$-variable polynomial.
\end{thm}
\begin{proof}
Let $Y_i$ be a Young diagram with at most two columns of type $(n-i,i)$, with $0\leq i\leq [n/2]$.
Let $p_i$ denote the corresponding orthogonal minimal central idempotent in $TL_n$, considered as an algebra over $\tZ$.
For example, $p_0$ is the Jones-Wenzl idempotent.

Since the full twist $\Delta^2$ is in the center of $B_n$, its image in any irreducible representation is a scalar.
The coefficient is a monomial that depends on the Young diagram, which can be computed using Lemma 9.3 \c{J} with $t=A^{-4}$.
Namely, the monomial $t^{k_i}$ for $Y_i$ as above has $k_i = i(n -i+1)$.
Therefore, the full twist can be represented in $TL_n$ as 
\[\rho(\Delta^2) = \sum_{i=0}^{[n/2]} t^{k_i} p_i \]
Since $p_i$ are orthogonal idempotents,
\[ \rho(\Delta^{2m}) = \left( \sum_{i=0}^{[n/2]} t^{k_i} p_i \right)^m = \sum_{i=0}^{[n/2]} t^{mk_i} p_i \]
If we express these idempotents in the $TL_n$ basis $\{h_j\}$, with $u_{ij}\in\tZ$, we can evaluate the bilinear pairing to derive the expression for $\langle L_m\rangle $:
\[ p_i = \sum_{j=1}^{C_n} u_{ij} h_j \]
\begin{equation}\label{Lm}
 \langle L_m\rangle  = \langle \Delta^{2m}, L''\rangle  = \sum_{i=0}^{[n/2]} t^{mk_i} \langle p_i, L''\rangle = \sum_{i=0}^{[n/2]} t^{mk_i} \sum_{j=1}^{C_n} u_{ij} \langle h_j, L''\rangle
\end{equation}
We define the rational function $P(t,x)$, which depends only on $L$ and $n$:
\begin{equation}\label{Ptx}
P(t,x)= \sum_{i=0}^{[n/2]} x^{k_i} \langle p_i, L''\rangle
\end{equation}
Therefore, 
\[ P(t,t^m) = \langle L_m\rangle  \qquad {\rm and} \qquad P(t,1) = \langle L\rangle . \]
The Jones polynomial equals the bracket up to a monomial depending on writhe, 
which does not affect the Mahler measure, so we obtain 
\[M(V_{L_m}(t)) = M(\langle L_m\rangle ) = M\left(P(t,t^m)\right)\]
We can now apply the special case of Lemma \ref{BL}:
\[ \lim_{m\to\infty} M(V_{L_m}(t))= \lim_{m\to\infty} M(P(t,t^m)) = M(P(t,x)) \]
By the proof of Theorem \ref{blocks}, $\tr(p_i)\in \frac{1}{\delta}\Z[\delta]$, so $(1+t) P(t,x)\in\Z[t^{\pm 1},x]$, which has the same Mahler measure as $P(t,x)$.
\end{proof}

$P(t,x)$ determines $V_{L_m}(t)$ for all $m$, and it is an interesting open question how it is related to the Jones polynomial of the link $L\cup U$.

In \c{Y}, Yokota used representation theory of the braid group to provide twisting formulas for the Jones polynomial.
Theorem \ref{thm1} follows from his Main Theorem by using Lemma \ref{BL} as in the proof above.
By introducing skein theory, though, we have simplified the argument, and extended it to colored Jones polynomials (see Theorem \ref{colJPthm}).  
Equation (\ref{Lm}) also determines a decomposition of the Jones polynomial into blocks after sufficiently many twists (see Theorem \ref{blocks}).

The proof of Theorem \ref{thm1} can be extended to produce a multivariable polynomial to which the Jones polynomial converges in Mahler measure under multiple twisting.
Given any link diagram $L$, construct $L_{m_1,\ldots,m_s}$ by surgeries: for $i=1,\ldots,s$, perform a $-1/m_i$ surgery on an unknot $U_i$ which encircles $n_i$ strands of $L$.
\begin{cor} \label{multitwist}
Let $\mathbf{m}=(1,m_1,\ldots,m_s)$. Let $L_\mathbf{m}$ be the multi-twisted link $L_{m_1,\ldots,m_s}$. 
The Mahler measure of the Jones polynomial of $L_\mathbf{m}$ converges as $\nu(\mathbf{m})\to\infty$ 
to the Mahler measure of an $(s+1)$-variable polynomial.
\end{cor}
\proof
We can suppose that the $U_i$ are far apart so these surgeries on $L$ are independent in the sense of the decomposition above.
In other words, if $L$ is originally given by $\bigcup_i \mathbf{1}_{n_i} \cup L''$, then the $i$th surgery replaces $\mathbf{1}_{n_i}$ by $\Delta^{2m_i}$.
Inductively, if $L_{m_0} = L''$, then for $k=1,\ldots,s$, we evaluate the bilinear pairing on $TL_{n_k}$:
\[ \langle L_{m_1,\ldots,m_k}\rangle  = \langle \Delta^{2m_k}, L_{m_1,\ldots,m_{k-1}}''\rangle  \]
Iterating the proof of Theorem \ref{thm1}, $P(t,x_1,\ldots,x_s)$ is constructed such that
\[ P(t,t^{m_1},\ldots,t^{m_s}) = \langle L_{m_1,\ldots,m_s}\rangle \]
We now directly apply Lemma \ref{BL},
$$ \lim_{\nu(\mathbf{m})\to\infty} M(V_{L_\mathbf{m}}(t)) = M\left(P(t,x_1,\ldots,x_s)\right) \eqno{\qed}$$

We now extend our results to the colored Jones polynomials.  Let $J_N(L; t)$ be the colored Jones polynomial of $L$, colored by the $N$-dimensional irreducible representation of $\SL_2(\C)$, with the normalization $J_2(L;t)=(t^{1/2}+t^{-1/2})V_L(t)$.  The colored Jones polynomials are weighted sums of Jones polynomials of cablings, and the following formula is given in \c{KM}.  Let $L^{(r)}$ be the $0$-framed $r$-cable of $L$; i.e., if $L$ is $0$-framed, then $L^{(r)}$ is the link obtained by replacing $L$ with $r$ parallel copies.
\begin{equation}\label{colJP}
J_{N+1}(L;t) = \sum^{[N/2]}_{j=0}(-1)^j\binom{N-j}{j}J_2({L^{(N-2j)}};t)
\end{equation}

\begin{thm}\label{colJPthm}
For fixed $N,$ and $L_m$ as above, the Mahler measure of $J_N({L_m};t)$ converges as $m\to\infty$, and similarly for multi-twisted links as $\nu(\mathbf{m})\to\infty$.
\end{thm}
\begin{proof}
Let $\Delta^2$ be the full twist $(\s_1\ldots\s_{k-1})^k$ in $B_k$.
The $0$-framed $r$-cable of $\Delta^2$ is the braid $(\s_1\ldots\s_{rk-1})^{rk}$, which is the full twist in $B_{rk}$.
Therefore, the operations twisting by full twists and $0$-framed cabling commute:
If $\tau_U^m(L)$ denotes $m$ full twists on the strands of $L$ encircled by an unknot $U$, then $\tau_U^m(L^{(r)})=(\tau_U^m(L))^{(r)}$.  
This is just a version of the ``belt trick'' (see, eg, I.2.3 \c{Prasolov-Sossinsky}).
So without ambiguity, let $L_m^{(r)}$ denote $\tau_U^m(L^{(r)})$.

Fix $r$ such that $1\leq r\leq N$.  
Following the proof of Theorem \ref{thm1}, we define $P_r(t,x)$ for the link $L^{(r)}$ such that $P_r(t,t^m) = V_{L_m^{(r)}}(t)$.
\[ \lim_{m\to\infty} M(V_{L_m^{(r)}}(t))= \lim_{m\to\infty} M(P_r(t,t^{m})) = M(P_r(t,x)) \]
We now apply (\ref{colJP}) to find $\widetilde{P}_N(t,t^m)$ for the $N$th colored Jones polynomial of $L_m$:
\begin{eqnarray*}
J_{N+1}({L_m};t) &=& \sum^{[N/2]}_{j=0}(-1)^j\binom{N-j}{j}(t^{\frac{1}{2}}+t^{-\frac{1}{2}})\,V_{L_m^{(N-2j)}}(t) \\
&=& \sum^{[N/2]}_{j=0}(-1)^j\binom{N-j}{j}(t^{\frac{1}{2}}+t^{-\frac{1}{2}})P_{N-2j}(t,t^m) \\
&=& (t^{\frac{1}{2}}+t^{-\frac{1}{2}})\,\widetilde{P}_{N+1}(t,t^m)
\end{eqnarray*}
Therefore,
$ \lim_{m\to\infty} M\left( J_N({L_m};t)\right)= \lim_{m\to\infty}M\left(\widetilde{P}_N(t,t^m)\right) =  M\left(\widetilde{P}_N(t,x)\right)$.
The proof for multi-twisted links follows just as for Corollary \ref{multitwist}.
\end{proof}


Many authors have considered the distribution of roots of Jones polynomials for various families of twisted knots and links (eg, \c{Z2,Z1,Z4,Z3}).  
The results above can be used to prove their experimental observation that the number of distinct roots approaches infinity, but for any $\epsilon > 0$, 
all but at most $N_{\epsilon}$ roots are within $\epsilon$ of the unit circle:

\begin{thm}\label{zeros}
Let $L_m$ be as in Theorem \ref{thm1}, or $L_\mathbf{m}$ as in Corollary \ref{multitwist}.  
Consider the family of polynomials, $V_{L_m}(t)$ or $J_N({L_m};t)$ for any fixed $N$, that vary as $m\to\infty$.
Let $\{\gamma^m_i\}$ be the set of distinct roots of any polynomial in this family.  
For any $\epsilon > 0$, there is a number $N_{\epsilon}$ such that 
\[ \#\{\gamma^m_i : \big| |\gamma^m_i| - 1 \big| \geq \epsilon \} < N_{\epsilon} \text{\; and\; } \lim\inf\nolimits_m\#\{\gamma^m_i \} \to \infty \text{\; as\; } m\to\infty, \]
and similarly for $L_\mathbf{m}$ as $\nu(\mathbf{m})\to\infty$.
\end{thm}
\begin{proof} 
We will prove the claim for $V_{L_m}(t)$, and the other cases follow similarly.  
By (\ref{Ptx}), $P(t,t^m) = \langle L_m\rangle$.  
So the $L^1$-norm of coefficients of $(1+t)V_{L_m}(t)$ equals that of the polynomial $(1+t)P(t,x)$, which is constant as $m\to\infty$.
This also follows from Theorem \ref{blocks}.
By Lemma 1 on p.187 of \c{Schinzel2}, any polynomial with a root $\gamma\neq 0$ of multiplicity $n$ has at least $n + 1$ non-zero coefficients.
Thus, for any integer polynomial $f$ with $f(0)\ne 0$, if the $L^1$-norm of its coefficients is bounded by $M$, and the number of its distinct roots is bounded by $k$, then $\deg(f)\leq k(M-1)$.
By (\ref{Lm}), the degree of $V_{L_m}(t)$ approaches infinity as $m\to\infty$.  
It follows that for any infinite sequence of polynomials in $\{V_{L_m}(t),\; m\geq 0\}$, there is a subsequence for which the number of distinct roots approaches infinity.

Let $\{\alpha^m_i\}$, for $1\leq i \leq a(m)$, be the roots counted with multiplicity of $V_{L_m}(t)$ outside the closed unit disc, $|\alpha^m_i|>1$.
Let $\{\beta^m_j\}$, for $1\leq j \leq b(m)$, be the roots counted with multiplicity inside the open unit disc, $|\beta^m_j|<1$.
Let
\[ A_{\epsilon}(m)= \#\{\alpha^m_i : |\alpha^m_i| - 1 \geq \epsilon \} \text{\; and\; } B_{\epsilon}(m)= \#\{\beta^m_j : 1 - |\beta^m_j| \geq \epsilon \} \]
By Theorem \ref{thm1} and Jensen's formula (\ref{Jensen}), $\lim\limits_{m\to\infty}\prod_{i=1}^{a(m)}|\alpha^m_i|$ exists.
Taking the mirror image $L_m^*$, $V_{L_m^*}(t)=V_{L_m}(t^{-1})$, or alternatively by the proof of Theorem \ref{thm1},
\[ \lim_{m\to\infty} M(V_{L_m}(t^{-1})) = M(P(t^{-1},t^{-m})) = M(P(t^{-1},x^{-1})) \]
so again by Jensen's formula (\ref{Jensen}), we have that $\lim\limits_{m\to\infty}\prod_{j=1}^{b(m)}|\beta^m_j|^{-1}$ exists.
\[\prod_{i=1}^{a(m)}|\alpha^m_i|\geq (1+\epsilon)^{A_{\epsilon}(m)} \text{\; and\; } \prod_{j=1}^{b(m)} |\beta^m_j|^{-1} \geq \left(\frac{1}{1-\epsilon}\right)^{B_{\epsilon}(m)} \]
Thus, there exist bounds $A_{\epsilon}(m)< A_{\epsilon}$ and $B_{\epsilon}(m)< B_{\epsilon}$.
Let $N_{\epsilon} = A_{\epsilon} + B_{\epsilon}$.
\end{proof}

\begin{example} \label{torus}
\rm
Let $T(m,n)$ be a torus knot, which is the closure of the braid $(\s_1\ldots\s_{n-1})^m$, with $m$ and $n$ relatively prime.
By 11.9 \c{J},
\[ V_{T(m,n)}(t)= \frac{t^{(n-1)(m-1)/2}}{(1-t^2)} (1-t^{m+1}-t^{n+1}+t^{n+m}) \]
Since the first factor has Mahler measure $1$, by Lemma \ref{BL},
\[ \lim_{m\to\infty} M(V_{T(m,n)}(t)) = M(1 - x \, t - t^{n+1} + x\, t^n) \]
Similarly,
\[ \lim_{m,n\to\infty} M(V_{T(m,n)}(t)) = M(1 - x \, t - y\, t + x\, y) \]
Since $M(V_{T(m,n)}(t))$ converges, by the proof of Theorem \ref{zeros}, the roots of $V_{T(m,n)}(t)$ approach the unit circle as $m+n\to\infty$.
This has been observed before (eg, \c{Z3}).
\end{example}

\begin{figure}
\begin{center}
\psfrag{dots}{\small{$\ldots$}}
\psfrag{ncross}{\small{n crossings}}
\psfrag{qk}{\small{$a_n$}}
 \includegraphics[height=1.5 in]{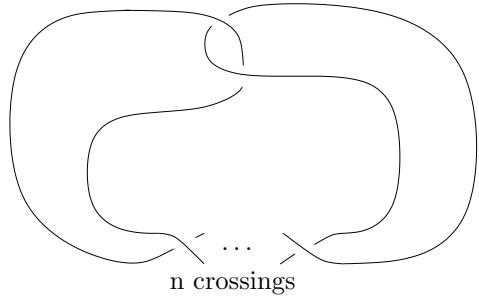}
\caption{The twist knot $K_n$}
\label{twistknot}
\end{center}
\end{figure} 

\begin{example} \label{twist}
\rm
Let $K_n$ be the twist knot shown in Figure \ref{twistknot}. 
Using Proposition \ref{2twist} below, the Jones polynomials up to multiplication by powers of $t$ and $\pm1$ are
$$V_{K_{2m}}(t)\doteq\frac{1-t^3+t^{2m+1}+t^{2m+3}}{1+t},\; \ V_{K_{2m+1}}(t)\doteq\frac{1-t^3-t^{2m+2}-t^{2m+4}}{1+t}$$
Since the denominator has Mahler meausre $1$, by Lemma \ref{BL},
\[ \lim_{m\to\infty} M(V_{K_{2m}}(t)) = M(1 - t^3+x t+x t^3)\] 
\[ \lim_{m\to\infty} M(V_{K_{2m+1}}(t)) = M(1 - t^3-x t^2 -x t^4) \]
Thus, as above, the roots of $V_{K_n}(t)$ approach the unit circle as $n\to\infty$.
\end{example}

\section{Twisting formulas}\label{twistsection}

Equation (\ref{Lm}) in the proof of Theorem \ref{thm1} actually provides an explicit structure for Jones polynomials after sufficiently many twists.
We observed this experimentally using a program written by Nathan Broaddus and Ilya Kofman.
This result extends to colored Jones polynomials.
If $V(t)=t^k\sum_{i=0}^s a_i t^i$, then $(a_0,\ldots,a_s)$ is the coefficient vector of $V(t)$.

\begin{thm} \label{blocks}
Suppose we twist $n$ strands of $L$ such that the Jones polynomial changes. 
For $m$ sufficiently large, the coefficient vector of $V_{L_m}(t)$ has $([n/2]+1)$ fixed possibly nontrivial blocks, one for each minimal central idempotent in $TL_n$, separated by blocks of zeros if $n$ is odd, or blocks of alternating constants, $\alpha,-\alpha,\alpha,-\alpha,\ldots$, if $n$ is even, which increase by constant length as $m$ increases.
\end{thm}
\begin{proof}
We rewrite (\ref{Lm}) by expressing $L''$ in the $TL_n$ basis, with $v_{ij}\in\ZA$,
\begin{equation}\label{Lm2}
\langle L_m\rangle  = \sum_{i=0}^{[n/2]} t^{mk_i} \langle p_i, L''\rangle = \sum_{i=0}^{[n/2]} t^{mk_i} \sum_{j=1}^{C_n} v_{ij} \langle p_i, h_j\rangle
\end{equation}
Using the Markov trace, $\langle p_i, h_j\rangle = \tr(p_i h_j)$.

The traces of minimal idempotents in $TL_n$ have been explicitly computed; for example, see Section 2.8 of \c{GHJ}.
However, to apply these trace formulas, we must renormalize according to our definition (\ref{TL}), and for $\tr(\rho(\beta))=\kb{\bar{\beta}}$, 
with the usual convention that for the unknot $\kb{\bigcirc}=1$.
Therefore, $\tr(\myone_n)=\delta^{n-1}$.
If $e_i'$ are the generators of $TL_n$ as in \c{J} and \c{GHJ}, then $e_i' = e_i/\delta$, and
\[ e_i' e_{i\pm 1}' e_i' = \tau e_i' \; \text{ where } \tau = \delta^{-2} \]
Let $P_k(x)$ be the polynomials defined by $P_0=P_1=1$ and for $k\geq 1$,
\[ P_{k+1}(x) = P_k(x) - x\, P_{k-1}(x) \]
According to Theorem 2.8.5 \c{GHJ}, with the renormalization, if $n$ is even, then $\tr(p_{n/2})=\tr(\myone_n) \tau^{n/2}$, and in other cases, for $i=0,1,\ldots,[n/2]$, 
\[ \tr(p_i) = \tr(\myone_n)\, \tau^i P_{n-2i}(\tau) \]
Observe that $P_{2k}$ and $P_{2k+1}$ have the same degree, but $\tr(\myone_{2k})$ and $\tr(\myone_{2k+1})$ will differ by a power of $\delta$.  
Therefore for all $i$, if $n$ is odd, $\tr(p_i)\in\Z[\delta]$, and if $n$ is even, $\tr(p_i)\in \frac{1}{\delta}\Z[\delta]$.
For example, the renormalized traces of the Jones-Wenzl idempotents $f_n$ in $TL_n$ are given by this formula as $p_0$ in $TL_n$:
\begin{eqnarray*}
\tr(f_2) = & \delta   \,  (1 - 1/\delta^2)       =  & (\delta^2 -1)/\delta \\
\tr(f_3) = & \delta^2 \,  (1 - 2/\delta^2)       =  &  (\delta^2 -2) \\
\tr(f_4) = & \delta^3 \,  (1 - 3/\delta^2 + 1/\delta^4) = & (\delta^4 - 3\delta^2 + 1)/\delta  \\
\end{eqnarray*}
To compute $\tr(p_i h_j)$, we recall another method of computing the Markov trace given in Section 5 of \c{J}, using weights associated to each Young diagram.  
Since $h_j$ are basis elements of $TL_n$, it follows that $\tr(p_i h_j) = \eta\, \tr(p_i)$, where $\eta\in\Z[\delta]$.  
In other words, multiplication by $h_j$ changes the weights according to the change of basis.  
This is also the idea of Section 13 of \c{J}, where the plat closure of a braid is considered; here, we consider all possible closures.  
Consequently, $\tr(p_i h_j)\in\Z[\delta]$ whenever $\tr(p_i)\in\Z[\delta]$.
For odd $n$, the result now follows from (\ref{Lm2}).

For even $n$, $\tr(p_i h_j)\in \frac{1}{\delta}\Z[\delta]$.  
All the idempotents add to $\myone$, so we can write the Jones-Wenzl idempotent 
$p_0= \myone - \sum_{i=1}^{[n/2]}p_i$.
\begin{eqnarray*}
\kb{L_m} & =& \sum_{i=0}^{[n/2]} t^{mk_i} \sum_{j=1}^{C_n} v_{ij} \kb{p_i, h_j} \\
&=& \sum_{j=1}^{C_n} v_{0j} \kb{\myone,h_j} - \sum_{i=1}^{[n/2]}(1-t^{m k_i})\sum_{j=1}^{C_n}v_{ij}\kb{p_i,h_j}\\
&=& q_0(t) - \sum_{i=1}^{[n/2]}\frac{1-t^{m k_i}}{1+t} q_i(t)\\
\end{eqnarray*}
where $q_0(t)=\displaystyle{ \sum_{j=1}^{C_n} v_{0j} \kb{\myone,h_j}}$ 
and $\displaystyle{q_i(t)=(1+t) \sum_{j=1}^{C_n}v_{ij}\kb{p_i,h_j}}$ for $1\leq i\leq [n/2]$.
Since $\langle p_i, h_j\rangle = \tr(p_i h_j)\in \frac{1}{\delta}\Z[\delta]$ and $1/\delta=-A^{-2}/(1+A^{-4})=-\sqrt{t}/(1+t)$, for all $i$, the $q_i(t)$ are Laurent polynomials with a possible $\sqrt{t}$ factor.

We consider the above summands separately. 
Since $k_i=i(n-i+1)$ and $n$ is even, $k_i$ is even for all $i$. Observe that 
\[ \frac{1-t^{2\ell}}{1+t}=\left( \frac{1-(t^2)^\ell}{1-t^2} \right)(1-t)=\left(\sum_{j=0}^{\ell-1}t^{2j}\right)(1-t)= \sum_{j=0}^{2\ell-1}(-1)^j t^j \]
If $q_i(t)$ has coefficient vector $(a_0,\ldots,a_s)$, for sufficiently large $m$, each summand looks like
\[ \left( \sum_{j=0}^{mk_i-1} (-1)^j t^{j} \right) q_i(t)= \bar{q}_i(t)+ q_i(-1)\left(\sum_{j=s}^{mk_i-1} (-1)^j t^{j}\right) + t^{mk_i}\tilde{q}_i(t) \]
where $\bar{q}_i(t)$ and $\tilde{q}_i(t)$ are polynomials of degree $s-1$ that depend on $q_i(t)$.
Hence, for sufficiently large $m$, the coefficient vector for each summand is  
\[(\bar{a}_0, \ldots, \bar{a}_{s-1},\, \alpha, -\alpha, \alpha, -\alpha, \ldots, \alpha,\, \tilde{a}_0, \ldots, \tilde{a}_{s-1} )\]
where all the coefficients depend on $q_i(t)$ and the outer blocks are fixed. 

Since $k_i=i(n-i+1)$, for sufficiently large $m$, the fixed blocks of the $i$-th summand do not interact with 
those of the $(i-1)$ summand.  For example,
\[ q_0(t)- \frac{1-t^{m k_1}}{1+t} q_1(t)-\frac{1-t^{m k_2}}{1+t} q_2(t) \]
has the following coefficient vector, with constants $r_0, r_1, r_2$:
\[ (a_1, \ldots, a_{r_0},\, \alpha, -\alpha, \ldots, \alpha, -\alpha,\, b_1, \ldots, b_{r_1},\, \beta, -\beta, \ldots \beta, -\beta,\, c_1, \ldots, c_{r_2} ) \]
By induction, we obtain the result for even $n$. \end{proof}

\begin{cor} \label{colblocks}
Suppose we twist $n$ strands of $L$ such that the Jones polynomial changes. 
If both $N$ and $n$ are fixed, then for $m$ sufficiently large, 
the coefficient vector of the colored Jones polynomial $J_N(L_m;t)$ has fixed blocks separated by 
blocks of zeros which increase by constant lengths as $m$ increases. 
\end{cor}
\begin{proof} Using the same notation as in the proof of Theorem \ref{colJPthm},
by the cabling formula (\ref{colJP}) and the ``belt trick,''
\[ J_{N+1}({L_m};t) = \sum^{[N/2]}_{j=0}(-1)^j\binom{N-j}{j}(t^{\frac{1}{2}}+t^{-\frac{1}{2}})\,V_{L_m^{(N-2j)}}(t) \]
By Theorem \ref{blocks}, the coefficient vector of $V_{L_m^{(N-2j)}}(t)$ has fixed blocks separated by alternating constants or zeros. 
The factor $(t^{\frac{1}{2}}+t^{-\frac{1}{2}})$ makes all the alternating constants to be zeros in each summand. 
In each summand, the number of strands of $L_m^{(N-2j)}$ being twisted is $n(N-2j)$.
Following the proof of Theorem \ref{blocks}, for each $j$ we expand $V_{L_m^{(N-2j)}}(t)$ using coefficients $k^j_i = i(n(N-2j)-i+1)$ in (\ref{Lm2}).
The result now follows from the sum over all $j$ of the corresponding equations (\ref{Lm2}). \end{proof}

We can obtain explicit formulas for $\kb{L_m}$ if we express the idempotents in terms of the basis of $TL_n$. 
Such formulas can be used to obtain the two-variable polynomial which appears in the limit of Theorem \ref{thm1}, as in Example \ref{twist}.
We now do this for $n=2$, primarily to illustrate Theorem \ref{blocks}, although this useful formula appears to be little known.

\begin{prop}\label{2twist}
For twists on 2 strands, if $L_m$ is obtained by adding $\Delta^{2m}$ at a crossing $c$ of $L$, then by splicing 
$c$ as in the Kauffman bracket skein relation,
\[ \kb{L_m}=A^{2m}\left(A\splice+\big(\sum_{i=0}^{2m}(-1)^iA^{-4i}\big)A^{-1} \kb{\asymp} \right) \]
\end{prop}
\begin{proof}  
Adding $\Delta^{2m}$ at a crossing $c$ is the same as adding $\Delta^{2m+1}$ to the link $L'$ obtained 
by splicing $L$ at $c$, such that $\kb{L',\myone}=\splice$.  Similarly, $\kb{L',e_1}=\kb{\asymp}$.

The basis of $TL_2$ is $\myone$ and $e_1$. 
There are two minimal central idempotents in $TL_2$: 
the Jones-Wenzl idempotent $p_0$, and $p_1=(\myone-p_0)$.
Now, $p_0=\myone-e_1/\delta$, and $p_1=e_1/\delta$.
Let $\Delta^2$ be the full right twist in $B_2$. 
If $\rho:B_2 \rightarrow TL_2$, then from the skein relation, $\rho(\Delta)=A p_0 -A^{-3}p_1$ and 
hence $\rho(\Delta^{2m+1})= A^{2m+1}(p_0-A^{-8m-4}p_1)$. Using the expression for $p_0$ and $p_1$ 
we get
\begin{eqnarray*}
\rho(\Delta^{2m+1})
&=&A^{2m+1}(\myone - \frac{e_1}{\delta} - A^{-8m-4}\frac{e_1}{\delta})\\
&=&A^{2m+1}\Big(\myone + A^{-2} \big(\frac{1+A^{-8m-4}}{1+A^{-4}}\big)e_1\Big)\\
&=& A^{2m}\Big(A\myone + \big(\sum_{i=0}^{2m}(-1)^iA^{-4i}\big)A^{-1}e_1\Big)
\end{eqnarray*}
Since $\kb{L_m}=\kb{L',\Delta^{2m+1}}$, the result now follows. \end{proof}

\begin{example}
  \rm Consider the closure of the $6$-braid, $\overline{12}3 4 3
  \overline{2} 1 \overline{2} 3 \overline{4 5 4} 3 \overline{2 4} 3
  5$, which is an $11$-crossing knot.  We perform $m$ full twists on
  the last $5$ and $6$ strands in Tables \ref{t5} and \ref{t6},
  respectively, and we give the coefficient vector, span, and Mahler
  measure of the resulting Jones polynomials.  Note there are three
  minimial central idempotents in $TL_5$, and four in $TL_6$.  The
  Mahler measure does not appear to converge to that of the
  corresponding 2-component links: For the 5-strand example,
  $M(V_{L\cup U}(t))\approx 7.998$, and for the 6-strand example,
  $M(V_{L\cup U}(t))\approx 12.393$. \end{example}

\begin{table}
\caption{Jones coefficients for $m$ twists on $5$ strands}\label{t5}

\bigskip
1  -3   8  -14  19  -23  23  -20  16  -9  4  -1 

\medskip
$m= 0,\quad \text{ span } = 11,\quad  M(V_{L_m}(t))=  4.198479$ 

\bigskip
1  -1   2  -1   4  -9  18 -28  35 -41  40 -35  26 -15   7  -2 

\medskip
$m= 1,\quad   \text{ span }= 15,\quad    M(V_{L_m}(t))=  5.785077$ 

\bigskip
1  -1   2  -1   2  -1   1   0   0   2  -8  16 -23  20 -11  -5  22 -36  44 -43  33 -19   8  -2 

\medskip
$m= 2,\quad \text{ span } = 23,\quad    M(V_{L_m}(t))=  8.267849$ 

\bigskip
1  -1   2  -1   2  -1   1   0   0   0   0   0   0   0   2  -8  16 -23  20 -12   0   7  -6  -1  14 -29  40 -42  33 -19   8  -2 

\medskip
$m= 3,\quad \text{ span } =31, \quad   M(V_{L_m}(t))=8.362212$ 

\bigskip
1  -1   2  -1   2  -1   1   0   0   0   0   0   0   0   0   0   0   0   0   2  -8  16 -23  20 -12   0   7  -7   4  -1   1  -5  15 -29  40 -42  33 -19   8  -2 

\medskip
$m= 4,\quad \text{ span } =39, \quad   M(V_{L_m}(t))=9.132926$ 

\bigskip
\framebox{1  -1   2  -1   2  -1   1}   0   0   0   0   0   0   0   0   0   0   0   0   0   0   0   0   0   \framebox{2  -8  16 -23  20 -12   0   7  -7   4  -1}   0   0   0   \framebox{1  -5  15 -29  40 -42  33 -19   8  -2} 

\medskip
$m= 5,\quad \text{ span } =47, \quad M(V_{L_m}(t))=8.568872$ 

\bigskip
\framebox{1  -1   2  -1   2  -1   1}   0   0   0   0   0   0   0   0   0   0   0   0   0   0   0   0   0   0   0   0   0   0   0   0   0   0   0   0   0   0   0   0   0   0   0   0   0   0   0   0   0   0   0   0   0   0   0   0   0   0   0   0   0   0   0   0   0   0   0   0   0   0   0   0   0   0   0   0   0   0   0   0   0   0   0   0   0   0   0   0   0   0   0   0   0   0   0 \framebox{2  -8  16 -23  20 -12   0   7  -7   4  -1}   0   0   0   0   0   0   0   0   0   0   0   0   0   0   0   0   0   0   0   0   0   0   0   0   0   0   0   0   0   0   0   0   0   0   0   0   0   0   0   0   0   0   0   0   0   \framebox{1  -5  15 -29  40 -42  33 -19   8  -2} 

\medskip
$m= 19,\quad \text{ span } =159,\quad   M(V_{L_m}(t))=8.589137$ 

\bigskip
\framebox{1  -1   2  -1   2  -1   1}   0   0   0   0   0   0   0   0   0   0   0   0   0   0   0   0   0   0   0   0   0   0   0   0   0   0   0   0   0   0   0   0   0   0   0   0   0   0   0   0   0   0   0   0   0   0   0   0   0   0   0   0   0   0   0   0   0   0   0   0   0   0   0   0   0   0   0   0   0   0   0   0   0   0   0   0   0   0   0   0   0   0   0   0   0   0   0   0   0   0   0   0   \framebox{2  -8  16 -23  20 -12   0   7  -7   4  -1}   0   0   0   0   0   0   0   0   0   0   0   0   0   0   0   0   0   0   0   0   0   0   0   0   0   0   0   0   0   0   0   0   0   0   0   0   0   0   0   0   0   0   0   0   0   0   0   0   \framebox{1  -5  15 -29  40 -42  33 -19   8  -2} 

\medskip

$m= 20,\quad \text{ span } =167,\quad\  M(V_{L_m}(t))=8.630147$
\end{table}

\begin{table}
\caption{Jones coefficients for $m$ twists on $6$ strands}\label{t6}

\bigskip
1   0   1   0   1   5 -14  23 -27  18  -4 -17  34 -46  49 -40  26 -11   2 

\medskip
$m= 1, \quad \text{ span } = 18,\quad   M(V_{L_m}(t))=  9.610867$ 

\bigskip
1   0   1   0   1   0   1  -1   1  -1   1   4 -14  22 -25  22 -27  37 -47  43 -21 -16  58 -98 127 -133 113 -74  36 -12   2 

\medskip
$m= 2, \quad \text{ span } = 30, \quad   M(V_{L_m}(t))=  15.131904$ 

\bigskip
1   0   1   0   1   0   1  -1   1  -1   1  -1   1  -1   1  -1   1   4 -14  22 -25  22 -27  36 -45  47 -45  45 -46  40 -20 -16  59 -105 150 -182 190 -167 123 -75  36 -12   2  

\medskip
$m= 3, \quad \text{ span } = 42,\quad    M(V_{L_m}(t))=  17.775295$ 

\bigskip
\framebox{1   0   1   0   1   0}   1  -1   1  -1   1  -1   1  -1   1  -1   1  -1   1  -1   1  -1   1  -1   1  -1   1  -1   1  -1   1  -1   1  -1   1  -1   1  -1   1  -1   1  -1   1  -1   1  -1   1  -1   1  -1   1  -1   1  -1   1  -1   1  -1   1  -1   1  -1   1  -1   1  -1   1  -1   1  -1   1  -1   1  -1   1  -1   1  -1   1  -1   1  -1   1  -1   1  -1   1  -1   1  -1   1  -1   1  -1   1  -1   1  -1   1  -1   1  -1   1  -1   1  -1   1  -1   1  -1   1  -1  \framebox{1 4 -14  22 -25  22 -27  36 -45  47 -45}  44 -44  44 -44  44 -44  44 -44  44 -44  44 -44  44 -44  44 -44  44 -44  44 -44  44 -44  44 -44  44 -44  44 -44  44 -44  44 -44  44 -44  44 -44  44 -44  44 -44  44 -44  44 -44  44 -44  44 -44  44 -44  44 -44  44 -44  44 -44  44 -44  44 -44  44 -44  44 -44  \framebox{45 -46  40 -20 -16  59 -105 151 -189 214 -223} 224 -224 224 -224 224 -224 224 -224 224 -224 224 -224 224 -224 224 -224 224 -224 224 -224 224 -224 224 -224 224 -224 224 -224 \framebox{223 -217 200 -168 123 -75  36 -12  2}

\medskip
$m= 19, \quad \text{ span } = 234,\quad    M(V_{L_m}(t))=  17.646099$ 

\bigskip
\framebox{1   0   1   0   1   0} 1  -1   1  -1   1  -1   1  -1   1  -1   1  -1   1  -1   1  -1   1  -1   1  -1   1  -1   1  -1   1  -1   1  -1   1  -1   1  -1   1  -1   1  -1   1  -1   1  -1   1  -1   1  -1   1  -1   1  -1   1  -1   1  -1   1  -1   1  -1   1  -1   1  -1   1  -1   1  -1   1  -1   1  -1   1  -1   1  -1   1  -1   1  -1   1  -1   1  -1   1  -1   1  -1   1  -1   1  -1   1  -1   1  -1   1  -1   1  -1   1  -1   1  -1   1  -1   1  -1   1  -1   1  -1   1  -1   1  -1  \framebox{1 4 -14  22 -25  22 -27  36 -45  47 -45}  44 -44  44 -44  44 -44  44 -44  44 -44  44 -44  44 -44  44 -44  44 -44  44 -44  44 -44  44 -44  44 -44  44 -44  44 -44  44 -44  44 -44  44 -44  44 -44  44 -44  44 -44  44 -44  44 -44  44 -44  44 -44  44 -44  44 -44  44 -44  44 -44  44 -44  44 -44  44 -44  44 -44  44 -44  \framebox{45 -46  40 -20 -16  59 -105 151 -189 214 -223} 224 -224 224 -224 224 -224 224 -224 224 -224 224 -224 224 -224 224 -224 224 -224 224 -224 224 -224 224 -224 224 -224 224 -224 224 -224 \framebox{223 -217 200 -168 123 -75  36 -12   2}

$m= 20, \quad \text{ span } = 246,\quad    M(V_{L_m}(t))=  17.622089$ 

\end{table}

\section{Pretzel Links}\label{pretzelsection}

\begin{figure}
\begin{center}
\psfrag{q1}{\small{$a_1$}}
\psfrag{q2}{\small{$a_2$}}
\psfrag{qk}{\small{$a_n$}}
 \includegraphics[height=2 in]{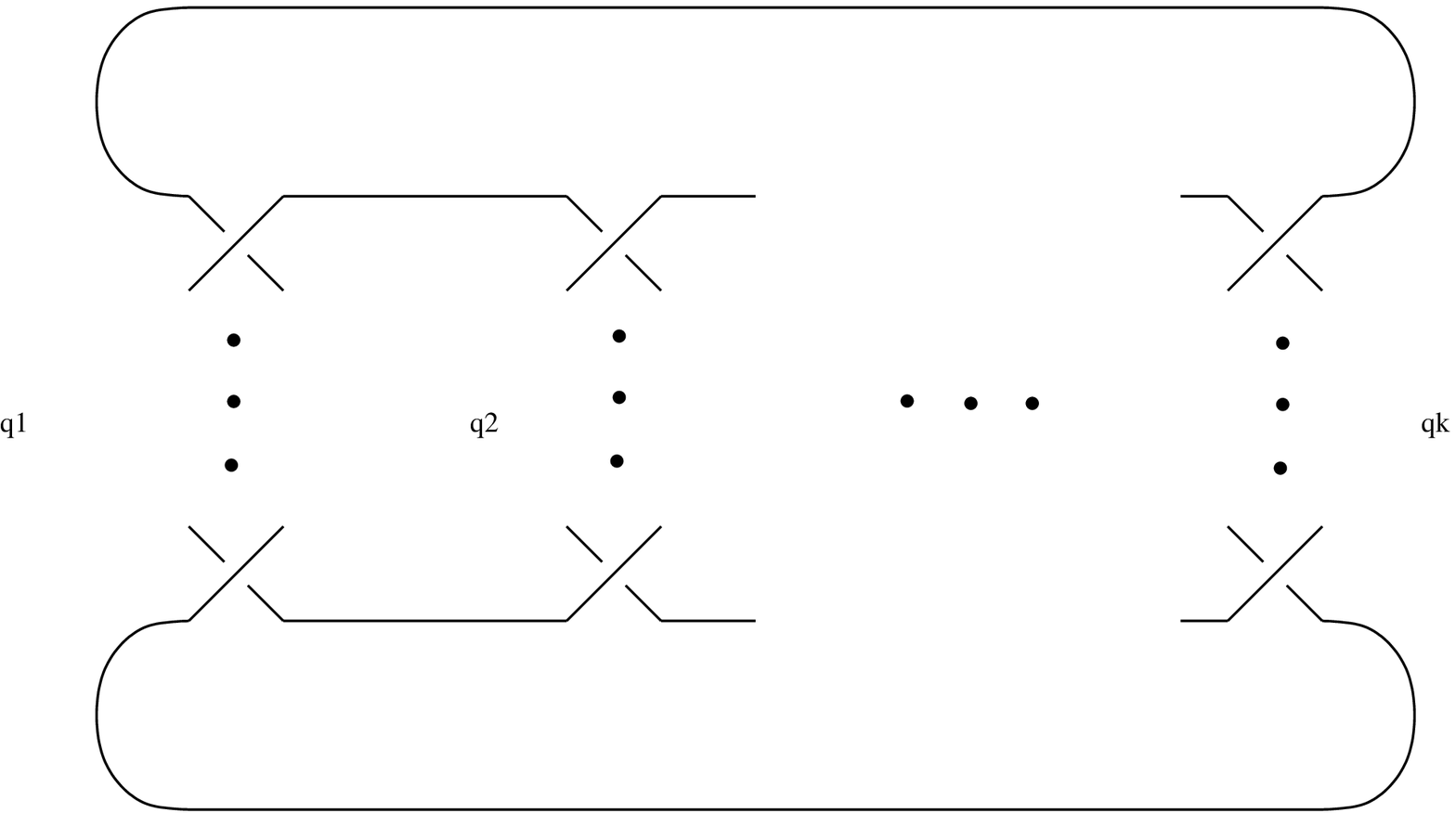}
\caption{The pretzel link $\mathcal{P}(a_1, \ldots, a_n)$}
\label{pretzelfig}
\end{center}
\end{figure}

In this section, we provide further evidence for the relationship between hyperbolic volume and Mahler measure of the Jones polynomial. 
Let $\mathcal{P}(a_1, \ldots, a_n)$ denote the the pretzel link as shown in Figure \ref{pretzelfig}. 
If $n$ is fixed, by Thurston's hyperbolic Dehn surgery theorem, the volume converges if all $a_i\to\infty$.
Similarly, by Corollary \ref{multitwist}, the Mahler measure of the Jones polynomial of the pretzel link $\mathcal{P}(a_1,\ldots,a_n)$ converges if all $a_i\to\infty$.

The main result of this section is that the Mahler measure of the Jones polynomial approaches infinity for $a_i=$ constant if $n\to\infty$, just as hyperbolic volume, according to Lackenby's lower bound in \c{Lackenby}, since the ``twist number'' is $n$ for $\mathcal{P}(a_1, \ldots, a_n)$.
In general, for an alternating hyperbolic link diagram, Lackenby gave lower and upper volume bounds using the twist number $T$ of the diagram.
Lackenby's lower bound for volume implies that as $T\to\infty$, Vol$(K_T)\to\infty$.
 Dasbach and Lin \c{DL} showed that $T$ is the sum of absolute values of the two Jones coefficients next to the extreme coefficients. 

Fix an integer $k \geq 1$.
Let $\mathcal{P}_n=\mathcal{P}(a_1, \ldots, a_n)$ where 
$a_1= \ldots =a_n=2k+1$, so $\mathcal{P}_n$ is a knot if $n$ is odd, and a link of two components if $n$ 
is even. 
 
\begin{thm}\label{pretzelthm}
$M(V_{\mathcal{P}_n}(t))\rightarrow \infty$ as $n \rightarrow \infty$.
\end{thm}
\begin{proof} Let $T$ denote the torus knot $T(2,2k+1)$. Let 
$T_n$ denote the connected sum of $T$ with itself $n$ times. Let 
$X=\displaystyle{ \sum_{i=0}^{2k}(-1)^iA^{-4i}}$. 
Let $Y=-A^4-1+X$.
Then 
\[ \kb{T}=A^{2k-1}\, Y \; \text{ and } \; \kb{T_n}=A^{(2k-1)n}\, Y^n\]
Consider $\mathcal{P}_n$ as the pretzel link $\mathcal{P}(2k+1, \ldots, 2k+1, 1)$ with $\Delta^{2k}$ added to the last crossing. 
Splicing this last crossing in two ways, we obtain $T_{n-1}$ and $\mathcal{P}_{n-1}$. By Proposition \ref{2twist},
$$\kb{\mathcal{P}_n}=A^{2k}\big( A \kb{T_{n-1}} + XA^{-1}\kb{\mathcal{P}_{n-1}} \big).$$
For a recursive forumula for any pretzel link, see \c{Landvoy}. 
Since $\mathcal{P}_2=T(2,4k+2)$,
\[\kb{\mathcal{P}_2}=A^{4k}\left( Y-A^{(-8k-4)}\, X \right)\]
By the recursion above and formulas for $\kb{T_n}$ and $\kb{\mathcal{P}_2}$,
\[\kb{\mathcal{P}_n}= A^{2nk-8k-n-2}\left(\frac{A^{(8k+4)}\, Y}{(-A^4-1)}\left( Y^{n-1}-X^{n-1} \right)\ -\ X^{n-1} \right) \]
Substitute $t=A^{-4}$. We now have 
\[ X=\sum_{i=0}^{2k}(-1)^it^{i}=\frac{1+t^{2k+1}}{1+t}\;  \text{ and } \; Y = -A^4-1+X=\frac{t^{2k+2}-t^2-t-1}{t(1+t)} \]
We substitute this in $\kb{\mathcal{P}_n}$. 
Up to multiplication by $\pm 1$ and powers of $t$, 
$$V_{\mathcal{P}_n}(t)\doteq \frac{(t^2+t+1)\big( t^{2k+2}+t \big)^n\ +\ 
t\big( t^{2k+2}-t^2-t-1\big)^n }{(1+t)^{n+1}}$$
Since the Mahler measure of the denominator is $1$,
\[ M\left(V_{\mathcal{P}_n}(t)\right)=M\left((t^2+t+1)\big( t^{2k+2}+t \big)^n\ +\ t \big( t^{2k+2}-t^2-t-1\big)^n \right).\]
By Lemma \ref{Boyd} below, $M(V_{\mathcal{P}_n}(t)) \rightarrow \infty$ as $n \rightarrow \infty$. \end{proof}

\begin{lem}\label{Boyd} 
For fixed integer $k \geq 1$,
$$M((t^2+t+1)\big( t^{2k+2}+t \big)^n\ +\ t \big( t^{2k+2}-t^2-t-1\big)^n ) \geq
c_k \ \rho^n $$ for some constants $\rho > 1$ and $c_k$.
\end{lem}
\begin{proof} The idea of the proof is due to David Boyd.  We define
\[a(t)= t^{2k+2}+t,\quad b(t)=t^{2k+2}-t^2-t-1,\quad c(t)=t^2+t+1,\quad
d(t)=t\]
\[p_n(t)=c(t)\,a(t)^n\ +\ d(t)\,b(t)^n\]
Let $I=[0,1]$.
Let $\beta\subset I$ be the union of intervals such that $|b(e^{2\pi i t})|
\leq |a(e^{2\pi i t})|$.  
To simplify notation, we omit the variable $e^{2\pi i t}$ from the polynomials below.
If $\alpha = b/|a|$,
\[ \log  \big| |c| - |\alpha|^n \big| \leq \log \frac{|p_n|}{\max\{|a|,|b|\}^n}  \leq \log \big( |c| + |\alpha|^n \big) \]
On $\beta,\, |\alpha|\leq 1$ and $|\alpha|= 1$ only at isolated points. 
Also, $|c| \pm |\alpha|^n=0$ only at isolated points, so these bounds are in $L^1(\beta)$.  
By the Dominated Convergence Theorem,
\[ \int_\beta \log\frac{|p_n|}{\max\{|a|,|b|\}^n} \rightarrow \int_\beta \log |c| \quad \text{ as } n\to \infty \]
Similarly for $I\setminus\beta$, so there exists a constant $C$ such that for $n$ sufficiently large,
\begin{equation} \label{C} 
 \int_I \log\frac{|p_n|}{\max\{|a|,|b|\}^n} > C
\end{equation}
Therefore,
\begin{eqnarray*}
\int_I \log|p_n| & > & n\,  \int_I \log(\max\{|a|,|b|\}) + C \\
\log M(p_n) & > & n\,  \log M(\max\{|a|,|b|\}) + C \\
M(p_n) & > & e^C\, \big(M(\max\{|a|,|b|\})\big)^n 
\end{eqnarray*}
Since $b(1) < 0 < b(2),\; b(t)$ has a root in the interval $(1,2)$, so $M(b)>1$.
Therefore, $M(\max\{|a|,|b|\})\geq \max\{M(a),M(b)\} > 1$.
\end{proof}

\begin{rmk}
\rm
We make some observations for knots with $\leq 16$ crossings using the Knotscape census. 
Only $17$ knots have Jones polynomials with $M(V_K(t))=1$, and they share only $7$ distinct Jones polynomials.
We list these in the table below using Knotscape notation; for example, $11n_{19}$ is the $19th$ non-alternating $11$-crossing knot in the census. 
Knots in the same row have the same Jones polynomial:

\leftskip25pt\begin{tabular}{lllll}
$4a_{1}$ & $ 11n_{19}$ & & & \\
$8a_{16}$ & $ 12n_{462}$ & $ 14n_{8212}$ & $ 16n_{509279}$ & \\
$9n_{4}$ & $ 16n_{207543}$  & & & \\
$12n_{562}$ & $ 12n_{821}$ & $ 13n_{1131}$ & $ 15n_{47216}$ & $ 16n_{23706}$ \\
$14a_{19115}$ & $ 16n_{992977}$  & & & \\ 
$14n_{26442}$  & & & & \\
$16a_{359344}$  & & & &
\end{tabular}

\leftskip0pt
Three pairs of these knots also have the same Alexander polynomial.
From the Knotscape pictures, it is easy to see that $16n_{23706}$ and $15n_{47216}$ are obtained from $11n_{19}$ by twisting in one and two places, respectively.
An interesting open question is how to construct more knots with $M(V_K(t))=1$.

Two knots in the census with the next smallest $M(V_K(t))$ are the following:

\leftskip25pt\begin{tabular}{ll}
$15n_{142389}$   &   $M(V_K(t))=  1.227786\ldots $  \\ 
$16a_{379769}$   &   $M(V_K(t))=  1.272818\ldots $   
\end{tabular}

\leftskip0ptSix knots have the next smallest value, 
$ M(V_K(t))=M(x^3-x-1)=1.324718\ldots $
Smyth showed that this polynomial has the smallest possible Mahler measure above 1 for non-reciprocal polynomials \cite{Smyth}.

For example, the torus knot $T(3,5)$ is one such knot, and it has the unusual property that it is also a pretzel knot, $\mathcal{P}(-2,3,5)$ (see Theorem 2.3.2 \c{Kawauchi}).
Torus knots are not hyperbolic, but with one more full twist on two strands, we obtain the twisted torus knot $T(3,5)_{2,1}$, which is also the pretzel knot $\mathcal{P}(-2,3,7)$, and this knot is the second simplest hyperbolic knot, with 3 tetrahedra.
$T(3,5)$ is $10_{124}$ in Rolfsen's table, the first non-alternating 10-crossing knot.  
Up to $10$ crossings, the next smallest $M(V_K(t))=1.360000$ belongs to $10_{125}$, which is hyperbolic with 6 tetrahedra.
$10_{125}$ is the pretzel knot $\mathcal{P}(2,-3,5)$, but apart from their Mahler measure, the Jones polynomials appear very different:

\begin{tabular}{llll}
 & \underline{degree} & \quad & \underline{coefficient vector} \\
$\mathcal{P}(-2,3,5)$ &  \ 4  \quad   10 & \quad & \ 1 \quad    0 \quad    1  \quad   0 \quad    0 \quad    0 \quad   -1  \\
$\mathcal{P}(2,-3,5)$ &  -4  \quad    \ 4 & \quad &  -1  \quad   1 \quad   -1 \quad    2 \quad   -1 \quad    2 \quad   -1 \quad    1 \quad   -1  
\end{tabular}\leftskip25pt

\end{rmk}

\Addresses\recd
\end{document}

%% file: agtout.tex
\def\ifplaintex{\expandafter\ifx\csname documentclass\endcsname\relax}

\def\gtp{{\mathsurround=0pt\it $\cal G\mskip-2mu$eometry \&\ 
$\cal T\!\!$opology $\cal P\!$ublications}}  

\def\recd{{\small Received:\qua\receiveddate\ifx\reviseddate\relax
\else\qquad Revised:\qua\reviseddate\fi\par}} 

\def\lognumber#1{\def\thelognumber{#1}}
\def\volumenumber#1{\def\thevolumenumber{#1}}
\def\volumeyear#1{\def\thevolumeyear{#1}}
\def\papernumber#1{\def\thepapernumber{#1}}
\def\pagenumbers#1#2{\def\startpage{#1}\def\finishpage{#2}}
\def\published#1{\def\publishdate{#1}}

\def\received#1{\def\receiveddate{#1}}
\def\revised#1{\def\reviseddate{#1}}
\def\accepted#1{\def\accepteddate{#1}}

\def\asciiaddress#1{\def\theasciiaddress{#1}}
\def\asciiemail#1{\def\theasciiemail{#1}}

\long\def\asciiabstract#1{\long\def\theasciiabstract{#1}}

\let\\\par\let\thelognumber\relax\let\thevolumenumber\relax
\let\thepapernumber\relax\let\thevolumeyear\relax\let\startpage\relax
\let\finishpage\relax\let\publishdate\relax\let\receiveddate\relax
\let\reviseddate\relax\let\accepteddate\relax\let\theasciititle\relax
\let\theasciiauthors\relax\let\theasciiaddress\relax
\let\theasciiabstract\relax

\let\theasciiemail\relax

\ifplaintex
\font\logobig=cmssbx10 scaled 3836
\font\logomed=cmssbx10 scaled 2557
\else
\font\logobig=cmssbx10 scaled 4200
\font\logomed=cmssbx10 scaled 2800
\fi

\long\def\makeagttitle{   
\count0=\startpage
\agt\hfill      
\hbox to 45truept{\vbox to 0pt{\vglue -13truept{\logomed A\kern -.37em{\logobig 
T}\kern -.38em G}\vss}\hss}
\break
{\small Volume \thevolumenumber\ (\thevolumeyear)
\startpage--\finishpage\nl
Published: \publishdate}

\vglue .25truein

{\parskip=0pt\leftskip 0pt plus
1fil\def\\{\par\smallskip}{\Large\bf\thetitle}\par\medskip} \vglue
0.05truein

{\parskip=0pt\leftskip 0pt plus 1fil\def\\{\par}{\sc\theauthors}
\par\medskip}%
 
\vglue 0.03truein 

{\small\leftskip 25truept\rightskip 25truept{\bf Abstract}\stdspace\theabstract

{\bf AMS Classification}\stdspace\theprimaryclass
\ifx\thesecondaryclass\relax\else; \thesecondaryclass\fi\par
{\bf Keywords}\stdspace \thekeywords\par}\vglue 7truept

}   

\ifplaintex
\font\phead=cmsl9 scaled 950
\font\pnum=cmbx10 scaled 913
\font\pfoot=cmsl9 scaled 950
\headline{\vbox to 0pt{\vskip -4.5mm\line{\small\phead\ifnum
\count0=\startpage ISSN 1472-2739 (on-line) 1472-2747 (printed)
\hfill {\pnum\folio}\else\ifodd\count0\def\\{ }%
\ifx\theshorttitle\relax\thetitle\else\theshorttitle\fi\hfill{\pnum\folio}
\else\def\\{ and }{\pnum\folio}\hfill\ifx\theshortauthors\relax\theauthors
\else\theshortauthors\fi\fi\fi}\vss}}
\footline{\vbox to 0pt{\vglue 0mm\line{\small\pfoot\ifnum\count0=\startpage
\copyright\ \gtp\hfill\else
\agt, Volume \thevolumenumber\ (\thevolumeyear)\hfill\fi}\vss}}
\else
\font=cmsl9 scaled 1050
\font\lnum=cmbx10 
\font=cmsl9 scaled 1050
\makeatletter
\def\@oddhead{{\small\lhead\ifnum\count0=\startpage ISSN 1472-2739 
(on-line) 1472-2747 (printed)\hfill {\lnum\number\count0}\else\ifodd\count0
\def\\{ }\ifx\theshorttitle\relax \thetitle \else\theshorttitle\fi\hfill
{\lnum\number\count0}\else\def\\{ and }{\lnum\number\count0}
\hfill\ifx\theshortauthors\relax 
\theauthors\else\theshortauthors\fi\fi\fi}}\def\@evenhead{\@oddhead}
\def\@oddfoot{\small\lfoot\ifnum\count0=\startpage\copyright\ \gtp\hfill\else
\agt, Volume \thevolumenumber\ (\thevolumeyear)\hfill\fi}
\def\@evenfoot{\@oddfoot}
\makeatother
\fi
\let\maketitlepage\makeagttitle
\let\makeshorttitle\maketitlepage
\let\maketitle\maketitlepage

\newwrite\gtoutfile
\long\gdef\makeheadfile{  
{\def\\{, }\def\s{ }
\immediate\openout\gtoutfile head.xxx
\immediate\write\gtoutfile{Proxy-for: \ifx\theasciiauthors\relax
\theauthors\else\theasciiauthors\fi\s<\ifx\theasciiemail\relax\theemail\else\theasciiemail\fi>}
\immediate\write\gtoutfile{\noexpand\\}
\immediate\write\gtoutfile{Authors: \ifx\theasciiauthors\relax
\theauthors\else\theasciiauthors\fi}
{\def\\{ }\immediate\write\gtoutfile{Title: \ifx\theasciititle\relax
\thetitle\else\theasciititle\fi}}
\immediate\write\gtoutfile{Subj-class: GT or SG, GR etc}
\immediate\write\gtoutfile{MSC-class: \theprimaryclass\ifx\thesecondaryclass\relax\else, \thesecondaryclass\fi}
\immediate\write\gtoutfile{Journal-ref: Algebr. Geom. Topol. \thevolumenumber\s
(\thevolumeyear) \startpage-\finishpage}
\immediate\write\gtoutfile{Comments: Published by Algebraic and
Geometric Topology at}
\immediate\write\gtoutfile{\s\s\s  http://www.maths.warwick.ac.uk/agt/AGTVol\thevolumenumber/agt-\thevolumenumber-\thepapernumber.abs.html}
\immediate\write\gtoutfile{\noexpand\\}
\immediate\write\gtoutfile{}
\ifx\theasciiabstract\relax
\immediate\write\gtoutfile{\theabstract}\else
\immediate\write\gtoutfile{\theasciiabstract}\fi
\immediate\write\gtoutfile{}
\immediate\write\gtoutfile{\noexpand\\}
\immediate\write\gtoutfile{}
\immediate\closeout\gtoutfile}}  

\def\maketitlepage{\makeagttitle\makeheadfile}
\let\makeshorttitle\maketitlepage
\let\maketitle\maketitlepage